\documentstyle[11pt]{article}

\newfont{\Bbb}{msbm10 scaled \magstep1}

\newcommand\bR{\hbox{\Bbb R}}
\newcommand\bZ{\hbox{\Bbb Z}}

\newcommand{\bA}{\hbox{\Bbb A}}
\newcommand{\bL}{\hbox{\Bbb L}}

\newfont{\es}{eusm10 scaled \magstep1}
\newfont{\ses}{eufm10 scaled \magstep1}
\newfont{\gt}{eufb10 scaled \magstep1}
\newfont{\goth}{eufb10 scaled \magstep2}

\newcommand{\gog}{\hbox{\gt g}}

\newcommand{\gn}{\hbox{\gt n}}

\newcommand{\goq}{\hbox{\gt q}}

\newcommand{\cw}{\hbox{\gt cw}}
\newcommand{\hog}{\hbox{\ses h}}
\newcommand{\dd}{\hbox{\ses d}}

\newcommand{\A}{\cal A}
\renewcommand{\L}{\cal L}
\newcommand{\B}{\cal B}

\newcommand{\I}{\cal I}

\newcommand{\W}{{\cal W}}

\def\ra{\rightarrow}
\def\hra{\hookrightarrow}
\def\be{\begin{equation}}
\def\ee{\end{equation}}
\newcommand\ve{\varepsilon}
\newcommand\tth{\tilde{\theta}}
\newcommand\tom{\tilde{\Omega}}
\newcommand\vfi{{\varphi}}

\newtheorem{theorem}{Theorem}[section]
\newtheorem{lemma}[theorem]{Lemma}
\newtheorem{proposition}[theorem]{Proposition}
\newtheorem{corollary}[theorem]{Corollary}
\newtheorem{remark}[theorem]{Remark}
\newtheorem{definition}[theorem]{Definition}
\newtheorem{ex}[theorem]{Example}

\begin{document}

\title{On a theorem of Henri Cartan concerning  the equivariant cohomology}

\author{Liviu I. Nicolaescu \\University of Notre Dame\\Notre Dame, IN 46556\\http://www.nd.edu/ $\tilde{}$ lnicolae/ }

\date{}

\maketitle

\begin{center}
{\bf Abstract}
\end{center}

{\it \small{We  provide a new look at an old result of Henri Cartan concerning the
cohomology of  infinitesimally free smooth Lie group actions}.}\footnote{{\bf 1991
Mathematics Subject Classification.} {\bf Primary} 57R91; {\bf Secondary} 58A10, 58A12}
\footnote{{\bf Key words and phrases}: Cartan and Weil models of equivariant cohomology,
algebraic connections, Chern-Weil transgression map.}

\bigskip

\addcontentsline{toc}{section}{Introduction}

\section*{Introduction}

 The equivariant cohomology of a smooth manifold  acted by a Lie group is a
 concept which crystallized in the works of  A. Borel, H. Cartan, C. Chevalley, H. Hopf, L. Koszul and A. Weil in the late forties and early fifties.

 The differential geometric approach to this subject was  brilliantly
 described by Henri Cartan in the beautiful survey \cite{Ca} which   continues to be  the first source for anyone  interested in learning  the basic facts of
  this theory.  Recently, it has been  the focus of intense research in connection with many problems in differential geometry, representation theory and  quantum field theory.

   A central result of  \cite{Ca}  is Cartan's theorem which
  states that if a compact  Lie group $G$ acts freely on a smooth manifold $M$
  then the $G$-equivariant  cohomology of $M$ (as defined by Cartan) is
  naturally isomorphic to the DeRham cohomology of the quotient.

  There are currently many proofs of this fact (e.g. \cite{Ca}, \cite{Du}, \cite{Ver}) but, in the author's view, they all suffer of the same  {\ae}sthetic  ``deficiency''.  They involve a  quite large amount of amazing  combinatorics whose origin is  somewhat obscure.  Moreover,  the resulting isomorphism is   extremely difficult to figure out {\em explicitely} at {\em cochain level}.

The main goal of this paper is to provide a new, {\em direct}  and
more transparent proof of the following  slight generalization of
Cartan's theorem.

    \bigskip

\noindent {\bf Theorem } {\em If $G$ is a compact Lie  group
acting on the
  smooth manifold $X$ and $N$ is a closed  normal subgroup of $G$ acting {\em
  freely} on  $X$ then the  $G$-equivariant cohomology of $X$ is {\em naturally}  isomorphic with the $G/N$-equivariant cohomology of $X/N$.}

    \bigskip

 We will actually establish a more general algebraic result (see Theorem
 \ref{th: Cartan}). Moreover, relying on a recent result of
Kalkman \cite{Ka1} (which provides a very explicit  isomorphism
between the
 Weil model and the Cartan model of equivariant cohomology) we will offer an {\em explicit} description of this isomorphism (along the lines of \cite{KV}). In the  course of the proof   we will  provide  yet another interpretation  for   {\em moment map}  and the equivariant characteristic classes described  for
 example in \cite{BGV}  or \cite{DV}.  The very simple functorial principle behind  our proof is explained in Remark \ref{rem: proof}.  While some of the computations involved may not look too eye pleasing, they are entirely routine and more importantly, their logical succesion is very natural.

There are two surprising aspects of this proof which make it so
attractive.  They can best be grasped  by looking at a special
example.   Suppose $P\ra B$ is a smooth principal $G$-bundle.
Cartan's theorem  then states that $H^*_G(P)$, the $G$-equivariant
cohomology of $P$, coincides with the  DeRham cohomology of the
base $B$. Naturally, one tries to construct cochain homotopy
equivalences between the complexes leading to the two
cohomologies.   A geometer might even attempt this  using purely
geometric operations on the smooth manifolds involved.  This
approach is doomed to fail.  The method we propose  is to embedd
the two complexes in  the same larger complex   consisting of
``ideal'' elements and then show that the two embeddings are
homotopic with this larger space.  The homotopies are described by
the Weil transgression between a  genuine $G$-connection on $P$
and a certain ``ideal'' connection which has only a formal meaning
!!!  The second surprise is the amazing effectiveness of this
method. Normally one expects that  by ``pushing'' the geometric
situation into an ``ideal'' abstract framework  the resulting
formul{\ae} will be more involved. To our surprize,   Kalkman's
isomorphism fits perfectly in such a framework.  A bonus of this
proof is that  the isomorphism $H^*_G(P) \cong H^*(B)$ can be {\em
explicitely} described at the  {\em cochain} level. More
precisely,   to   any Cartan representative of an element in
$H^*_G(P)$ we associate a closed form on  $B$. This correspondence
descends to an {\em isomorphism}  between the two cohomologies.
This map is obtained naturally, as a by-product of our
computations.

     The paper is divided into five section of which four are devoted  to surveying the basic ``players'' in Cartan's approach to equivariant theory. In the first part we introduce the notion of  {\em operation} which  captures the essential features of the DeRham algebra  of a smooth manifold with a Lie group action.

     In the first section we introduce the main object of study, that of {\em operation}. It formalizes the algebra of exterior forms on a smooth manifold equipped with a Lie group action. Section 2 introduces  the Cartan and Weil models of equivariant cohomology
 while  section 3 describes  the Kalkman isomorphism between them. In section 4 we review the basics of the Weil transgression trick in the framework of {\em operations}.  In the final section  we prove Cartan's theorem and  discuss a few consequences.

     In the sequel $G$ will always denote a compact, connected Lie group.  The prefix ``s''   will refer to ``super'' (i.e. ${\bZ}_2$-graded) objects  as  in
\cite{BGV}. The bracket $[\,\cdot \, , \, \cdot \,]_s$ denotes the
super-commutator  in a super-algebra. Also, we will use Einstein's
summation  convention (unless otherwise indicated).

\bigskip


\bigskip

\section{Operations}
\setcounter{equation}{0}

 As in \cite{DV} we will consider Frechet algebras. These are associative
 ${\bR}$-algebras such that their algebraic operations are continuous with
 respect to a Frechet topology. The standard example of Frechet algebra  is that of  the algebra of smooth functions on a
 smooth manifolds.

 In this section we want to introduce the algebraic counterpart of the
 geometric  notion of smooth manifold acted on by a Lie group.  This object
 appears in literature with various names. We have chosen the terminology  of
 \cite{GHV} which  stays closer to the original motivation.

 \begin{definition}{\rm An {\em operation} consists of the following.

 \noindent (a) A s-commutative ${\bZ}$-graded  Frechet algebra
 \[
 {\A}=\oplus_{n\in {\bZ}}{\A}^n
 \]
 such that $a\cdot b =(-1)^{|a|\cdot|b|}b\cdot a$ for any two homogeneous
 elements. (${\A}$ is naturally a s-algebra  by ${\A}={\A}^{even}\oplus {\A}^{odd}$).

 \noindent (b) A {\em continuous  odd} derivation
 \[
 d:{\A}^{*}\ra {\A}^{*+1}
 \]
 such that $d^2=0$.

 \noindent (c) A {\em smooth} action of the  Lie group $G$ on ${\A}$  via
 algebra automorphisms commuting with $d$. We denote by $\L$ the derivative of
 this action  at ${\bf 1}\in G$. Thus ${\L}$ defines a representation of
 ${\gog}$ into the Lie algebra of even derivations of $\A$.  (${\L}_X$  $(X\in
 {\gog})$ is the Lie derivative along the automorphisms $\exp(tX)$ of $\A$.)
 Note that
 \[
 [{\L}_X ,d]_s={\L}_Xd -d{\L}_X =0.
 \]

 \noindent (d) A continuous $G$-equivariant linear map $\I$ from ${\gog}$ (called {\em
 contraction}) to the space of {\sl
 odd} derivations of $\A$ such that $\forall \, X, Y\in {\gog}$

 (d1) ${\I}_X({\A}^n)\subset {\A}^{n-1},\;\;\forall n.$

 (d2) $[{\I}_X,{\I}_Y]_s={\I}_X {\I}_Y+{\I}_Y{\I}_X= 0$.

 (d3) $ [{\L}_X, {\I}_Y]= {\I}_{[X,Y]}$.

 (d4) (Cartan formula)} $[{\I}_X, d]_s={\I}_X d+d{\I}_X= {\L}_X$.
 \end{definition}

 \begin{remark}{\rm The above contraction ${\I}$ extends to an algebra morphism
 ${\I}:\Lambda {\gog} \ra {\rm End}\, ({\A})$ thus defining a $\Lambda
 {\gog}$-module structure on ${\A}$.}
 \end{remark}

 \begin{ex}{\rm The {\em right}  action of a Lie group $G$ on the smooth manifold $M$
 defines a structure of operation on $\Omega^*(M)$. For each $X\in {\gog}$ we
  will denote  by ${\L}_X$ the Lie derivative along the flow $m\mapsto m\cdot \exp(tX)$ while ${\I}_X$ denotes the contraction along $X^{\#}$-the infinitesimal generator of the above
 flow.}
 \end{ex}

 \begin{ex} {\bf (The Weil algebra)} {\rm  Consider the Lie group $G$ and  set
 \[
 W_G=S {\gog}^*\otimes \Lambda {\gog}^*
 \]
 where $\Lambda$ and  $S$ denote the exterior  and respectively
 the symmetric algebra.  Topologize $W_G$ in the obvious fashion (as a space of
 polynomials) and  equip it
 with the ${\bZ}$-grading
 \[
 \deg (\Lambda^p {\gog}^*) =p,\;\;\deg (S^q {\gog}^*) =2q.
 \]
 Denote by $\hog$ the obvious isomorphism
 \be
 {\hog}:\Lambda^1{\gog}^* \ra S^1 {\gog}^*.
 \label{eq: augment}
 \ee
 The usual derivations $d$, $L_X$ and $i_X$ on $\Omega^*(G)$ have an algebraic counterpart on  $\Lambda {\gog}^* $ which we denote by ${\bf d}$, ${\bf L}_X$ and $\imath_X$., It is convenient to describe these operations in ``local coordinates''. Choose a basis $(e_i)$ of ${\gog}$ and denote by $(\theta^i)$ the dual basis  of ${\gog}^*$. Set $\Omega^i={\hog}(\theta^i)\in S^1 {\gog}^*$. Denote by ${\bf L}_i$ and $\imath_j$ the Lie derivative along
 $e_i$ and respectively the contraction by $e_j$. Then
 \[
 {\bf L}_i \theta^j = -C^j_{ik}\theta^k\;\;{\rm
 and}\;\;\imath_j\theta^k=\delta^k_j
 \]
 where $C^i_{jk}$ denote the structural constants of the Lie algebra ${\gog}$
 \[
 [e_j, e_k]=C^i_{jk}e_i.
 \]
 To  describe {\bf d} in local coordinates we  use the formula
 \[
 ({\bf d} \omega)= -\omega([X,Y]),\;\;\omega\in {\gog}^*,\; X,Y \in {\gog}.
 \]
 This yields
 \[
 {\bf d}\theta^i = -\frac{1}{2}C^i_{jk}\theta^j \theta^k.
 \]
 This implies immediately Koszul's formula
 \be
 {\bf d}= \frac{1}{2} \mu(\theta^k ){\bf L}_k = \frac{1}{2} \theta^k {\bf L}_k
 \label{eq: koszul}
 \ee
 where for any algebra $A$ and any $a\in A$ we denote by $\mu(a)$ the left
 multiplication by  $a$.  For simplicity we will very often omit the $\mu$
 symbol when  there is very little room for confusion.

 The action of $G$ on $W_G$ induced by the coadjoint action  defines a Lie
derivative ${\bf L}$ extending the Lie derivative on $\Lambda
{\gog}^*$ according to the prescription
\[
{\bf L}_X {\hog}(\omega) ={\hog}({\bf L}_X \omega).
\]

 We can extend ${\bf d}$ to an odd derivation $d_W$ of $W_G$ uniquely determined by its action on the generators
 \be
 d_W\theta^i ={\bf d}\theta^i +{\hog}\theta^i
 =-\frac{1}{2}C^i_{jk}\theta^j\theta^k +\Omega^i = (\frac{1}{2}\theta^j{\bf L}_j +\Omega^j\imath_j)\theta^i
\label{eq: we1}
 \ee
 \be
 d_W\Omega^i=-C^i_{jk}\theta^j\Omega^k = \left(\theta^j {\bf L}_j \right)\Omega^i.
\label{eq: we2}
 \ee
Extend the contraction $\imath_X$ to $W_G$ by imposing $\imath_X
{\hog}(\omega)=0$.

We leave the reader  verify that  these three derivations on $W_G$
do indeed define a structure of operation.  $\Box$}
\end{ex}

Given two operations $({\A}_i, d_i, {\L}^i,  {\I}^i)_{i=1,2}$ we
can define a structure of operation on their (Grothendiek)
topological tensor product ${\A}_1\otimes {\A}_2$ (described e.g.
in \cite{Tre}) according to the rules $({\ve}_1$ is the grading
operator of the s-algebra ${\A}_1$)
\[
d =d_1 \otimes {\bf 1} +{\ve}_1\otimes d_2
\]
\[
{\L}={\L}^1\otimes {\bf 1} +{\bf 1}\otimes {\L}^2
\]
and
\[
{\I}={\I}^1 \otimes {\bf 1} +{\ve}_1\otimes {\I}^2.
\]

Given an operation $({\A}, d, {\L}, {\I})$ we can define three
subalgebras
\[
{\A}_{inv}=\ker {\L} =\{a \in {\A} \; ;\; {\L}_X a=0 \;\;\forall
X\in {\gog} \}
\]
\[
{\A}_{hor}=\ker {\I}\;\;{\rm and}\;\;{\A}_{bas}= {\A}_{inv}\cap
{\A}_{hor}.
\]
Since $[d, {\L}]_s =0$ we deduce $d{\A}_{inv}\subset {\A}_{inv}$.
Moreover, Cartan formula implies  $d{\A}_{bas}\subset {\A}_{bas}$.
Thus we can define the cohomology groups
\[
H^*_{inv}({\A})=H^*({\A}_{inv}, d)\;\;{\rm and}\;\;
H^*_{bas}({\A})=H^*({\A}_{bas},d).
\]
\begin{ex}{\rm Consider a smooth principal $G$-bundle $G{\hra}P \ra B$. The right action of $G$ on $P$ induces a structure of operation on $\Omega^*(P)$. The basic subalgebra of this  operation is then naturally isomorphic to $\Omega^*(B)$.}
\end{ex}

\section{The Cartan-Weil descriptions}
\setcounter{equation}{0}

We will work in the more general setting of  operations. We will
define two notions of  equivariant cohomology and in following
subsection we will show they coincide.

Consider a $G$-operation $({\A},d, {\L}, {\I})$.

\bigskip

\noindent {\bf The Weil description}  We define  Weil's
equivariant cohomology of  ${\A}$ by
\[
WH^*_G({\A}) \stackrel{def}{=} H^*_{bas}(W_G\otimes {\A})
\]
where $W_G$ denotes the Weil algebra introduced in the previous
subsection.

\bigskip

\noindent {\bf The Cartan description} Consider the algebra
\[
{\B}=S{\gog}^* \otimes {\A}.
\]
$S{\gog}^*$ is graded as usually by
\[
\deg S^p{\gog}^* =2p.
\]
$G$ acts smoothly on ${\B}$ and so we can form the  subalgebra of
invariant elements
\[
{\B}_{inv}=(S{\gog}^*\otimes {\A})^G.
\]
Now define the operator
\be
{\dd}={\bf 1}\otimes  d - \sum_i\mu(\Omega^k)\otimes {\I}_k={\bf
1}\otimes d- \Omega^k\otimes {\I}_k. \label{eq: dcar} \ee If we
regard $\omega \in S{\gog}^*\otimes {\A}$ as a polynomial map
${\gog}\ra {\A}$ then ${\dd}\omega$ is the polynomial map
\[
X\mapsto d(\omega(X))-{\I}_X(\omega(X)).
\]
 ${\dd}$ satisfies the following conditions (see \cite{Ca})
\be
{\dd}^2=-\Omega^k\otimes {\L}_k. \label{eq: dcar1} \ee
\be
[{\dd}, {\bf L}\otimes {\bf 1}+{\bf 1}\otimes {\L}]=0. \label{eq:
dcar2} \ee

The equality (\ref{eq: dcar2}) shows that  ${\B}_{inv}$ is ${\dd}$
invariant.  Moreover, on this  subalgebra  ${\dd}^2 =0$.   Indeed,
on this subalgebra  we have ${\bf L}\otimes {\bf 1} = - {\bf
1}\otimes {\L}$ so that by (\ref{eq: dcar1}) we have
\[
{\dd}^2= -\Omega^i\otimes {\L}_i = \Omega^i{\bf L}_i \otimes {\bf
1}.
\]
Now it is not  difficult to see that $ \Omega^i{\bf L}_i \equiv 0$
on $S{\gog}^*$ due
 to the skew symmetry of the structural constants. Thus $({\B}_{inv}, {\dd})$ is a cochain complex  and  we define the Cartan equivariant cohomology of  ${\A}$ by
\[
CH^*_G({\A}) \stackrel{def}{=}H^*({\B}_{inv}, {\dd}).
\]
When ${\A}$ is the algebra of differential forms on a smooth
manifold $M$ on which $G$ acts smoothly  we will use the notations
$WH^*_G(M)$ and $CH^*_G(M)$ to denote the corresponding
equivariant cohomologies.

\section{Weil model $\Longleftrightarrow $ Cartan model}
\setcounter{equation}{0}
 Consider a $G$-operation $({\A}, d, {\L}, {\I})$.   The main result of this
 subsection is the following.

 \begin{theorem}{\rm  There exists a} natural {\rm isomorphism}
 \[
 WH^*_G({\A})\cong CH^*_G({\A}).
 \]
\label{th: CarWe}
 \end{theorem}

We briefly describe the proof in \cite{Ka1}. For a different but
related approach we refer to \cite{MQ}.

 Consider the algebra ${\B}= W_G\otimes {\A}$.  It has a tensor product structure of $G$-operation with structural derivations    $D$, $L$ and respectively $I$.   For each $U^j_i=\theta^j \otimes e_i \in {\gog}^*\otimes {\gog}$ define  the following operators (on
 ${\B}$)
\[
{\bA}^j_i=\theta^j\otimes {\I}_i
\]
\[
{\bL}^j_i=\theta^j\otimes {\L}_i -\Omega^j\otimes {\I}_i
\]
In general for any $T=t^i_jU^j_i \in {\gog}^*\otimes {\gog}$ set
\[
{\bA}_T=t^i_j{\bA}^j_i,\;\;{\bL}_T= t^i_j{\bL}^j_i\;\;{\rm
and}\;\; D_T=D+{\bL}_T.
\]
Note that $({\bA}^j_i)^2=0$, $\forall i,j$ and moreover
${\bA}_T{\bA}_S={\bA}_S{\bA}_T$, $\forall  S, T$. Thus
$\exp({\bA}_T)$ is well defined and invertible.  A simple
computation shows that for all $i, j$  the operator
$\exp({\bA}^j_i)={\bf 1}+{\bA}^j_i$ is an {\em algebra
automorphism} of ${\B}$  so that $\exp({\bA}_T)$ is an
automorphism  of ${\B}$ for all $T\in {\gog}^*\otimes {\gog}$.

The key step in the proof of Theorem \ref{th: CarWe} is contained
in the following result.

\begin{lemma}{\rm For any $T\in {\gog}^*\otimes {\gog}$ we have}
\[
\exp({\bA}_T)D\exp(-{\bA}_T )=D_T.
\]
\label{lemma: brst}
\end{lemma}

\noindent{\bf Proof of the lemma}  An elementary computation
shows that for any $i, j,k,\ell$ we have the following
``differential equations''
\be
[D, \exp({\bA}^i_j)]={\bL}^i_j \exp({\bA}^i_j) \label{eq: brst1}
\ee
\be
[{\bL}^i_j , \exp({\bA^k_l})]=0. \label{eq: brst2} \ee The
equality  (\ref{eq: brst1}) can be rephrased as
\be
\exp({\bA}^i_j)D\exp(-{\bA}^i_j)=D_{U^i_j}:=D+{\bL}^i_j \label{eq:
brst3} \ee while (\ref{eq: brst2}) is equivalent to
\be
\exp({\bA}^k_\ell){\bL}^i_j\exp(-{\bA}^k_\ell)={\bL}^i_j.
\label{eq: brst4} \ee Using (\ref{eq: brst4}) in (\ref{eq: brst3})
we deduce $({\bA}={\bA}^i_j+{\bA}^k_\ell)$
\[
\exp({\bA})D\exp(-{\bA}) =D_{\bA}.
\]
Lemma \ref{lemma: brst} now follows by iterating the above
procedure.  $\Box$

\bigskip

We now have a whole family of  $G$-structures  on ${\B}$
parameterized  by ${\gog}^*\otimes {\gog}$,
\[
{\B}_T:=({\B}, D_T, L_T=\exp({\bA}_T)L\exp(-{\bA}_T),
I_T=\exp({\bA}_T)I\exp(-{\bA}_T) \, ),\;\;T\in
{\gog}^*\otimes{\gog}.
\]
Moreover  an elementary computation shows that $L_T\equiv  L_0$.
All these structures are isomorphic with the canonical tensor
product structure and in particular
\[
H^*_{bas}({\B}_T)\cong H^*_{bas}({\B}_0)\cong  WH^*({\A}).
\]
An interesting special case arises when $T={\bf
id}=\theta^i\otimes e_i$. In this case the derivation $D_{\bf id}$
is known as the {\em BRST} ($=$  Bechi-Rouet-Stora-Tyupin)
operator and it arises in the  quantization of classical gauge
theories.

We leave the reader check that
\[
({\B}_{\bf id})_{hor} = \ker I_{\bf id} =S{\gog}^* \otimes {\A}.
\]
Hence
\[
({\B}_{\bf id})_{bas}= (S{\gog}^*\otimes {\A})^G
\]
and   it is not difficult to see that on this subalgebra $D_{\bf
id}={\dd}$. Thus
\[
H^*_{bas}({\B}_{\bf id}) \cong CH^*_G({\A}).
\]
Theorem \ref{th: CarWe} is proved. $\Box$

\begin{remark}{\rm Kalkman's isomorphism
\be
\phi =\phi_G=\exp({\bA}_{\bf id}):{\cal B}_0  \ra {\cal B}_{\bf
id} \label{eq: kalk} \ee has a particularly nice form when
restricted to $({\cal B}_0)_{bas}$. It is uniquely determined by
the correspondences $\Omega^i \mapsto \Omega^i$ and $\theta^j
\mapsto 0$.} \label{remark: brst}
\end{remark}

\section{Algebraic connections}
\setcounter{equation}{0}

Among the possible actions of a Lie group on a smooth manifold the
free ones play a special role. Consider for example the  case of a
smooth  principal  $G$-bundle  $G{\hra}P\ra B$. Such actions admit
connections. Recall (see \cite{KN}) that a connection  on $P$ is
an equivariant splitting
\[
TP\cong {\cal V}P\oplus {\cal H}P
\]
where ${\cal V}P$ is the bundle spanned by the infinitesimal
generators of  the $G$ actions.   In fact, for any $p\in P$ the
correspondence
\[
{\gog}\ni X \mapsto X^\#_p
\]
identifies the fiber ${\cal V}_pP$ with ${\gog}$.

Alternatively, a splitting as above can be defined  by a vertical
projector i.e. a ${\gog}$-valued 1-form $\Theta \in {\gog}\otimes
\Omega^1(P)$  which is $G$-invariant and satisfies
\be
i_{X^\#}\Theta =X\;\;\forall X\in {\gog}. \label{eq: proj} \ee We
can regard this connection as a linear map
\[
\Theta:{\gog}^* \ni \theta^i\mapsto {\tth}^i
\]
so that
\[
\Theta =e_i\otimes {\tth}^i.
\]
The condition (\ref{eq: proj}) reads
\be
{\I}_j{\tth}^k=\delta^k_j \label{eq: conn1} \ee or equivalently,
\be
{\I}_X\Theta =X,\;\;\;\forall X\in {\gog}. \label{eq: conn4} \ee
The invariance implies
\[
{\bf L}_ke_i\otimes {\tth}^i +e_i \otimes {\L}_k{\tth}^i=0
\]
i.e.
\be
{\L}_k{\tth}^i=-C^i_{kj}{\tth}^j. \label{eq: conn2} \ee or
equivalently,
\be
\Theta {\bf L}_X ={\L}_X \Theta,\;\;\;\forall X\in {\gog}.
\label{eq: conn3} \ee
 The conditions  (\ref{eq: conn1})-(\ref{eq: conn3}) are formulated  using a language which involves only the structure of operation.   Thus we can define an abstract  notion  of {\em algebraic connection} on any $G$-operation
 ${\cal A}$ as a $G$-equivariant linear map ${\tth}: {\gog}^* \ra {\cal A}$ satisfying
 (\ref{eq: conn1})-(\ref{eq: conn3}).

\begin{ex}{\rm The inclusion ${\gog}^* {\hra}W_G$ defines an algebraic connection.}
\end{ex}

Consider now a $G$-operation  $({\A}, d, {\L}, {\I})$ equipped
with a connection
\[
\Theta: \theta^i \mapsto {\tth}^i.
\]
Define
\[
{\tom}^i= d{\tth}^i-\frac{1}{2}C^i_{jk}{\tth}^j{\tth}^k .
\]
The form
\[
{\tom}=e_i\otimes {\tom}^i \in {\gog}\otimes {\A}^2
\]
is independent of the basis $(e_i)$ and it is called {\em the
curvature of the connection}. An elementary computation shows that
\[
{\I}_k{\tom}^i= 0\;\;\;\forall i,k
\]
i.e. ${\tom}^i\in {\A}^2_{hor}$, $\forall i$.

The main algebraic implications of the existence of a connection
derive from the following decomposition result.

\begin{proposition}{\rm The connection induced map
\[
\Lambda{\gog}^*\otimes {\A}_{hor} \ra {\A},\;\;\; \theta^A\otimes
\omega \mapsto {\tth}^A\omega
\]
($\omega \in {\A}_{hor}$, $A$ is an ordered multi-index $(a_1,
a_2, \ldots)$ and $\theta^A=\theta^{a_1}\wedge \theta^{a_2}\wedge
\cdots $) is an isomorphism of graded algebras.}
\end{proposition}

\noindent{\bf Idea of proof}\hspace{.3cm} The map is clearly
injective. The surjectivity follows from the following simple
observation
\[
\forall \omega\in {\A},\;\; \omega -{\tth}^k {\I}_k\omega \in \ker
{\I}_k \;\;{\rm (no\;\;summation)}.\;\;\;\Box
\]
Thus we can {\em uniquely} represent any element $\omega \in {\A}$
as a polynomial
\[
\omega={\tth}^A\omega_A
\]
where in the above sum $A$ runs through all ordered multi-indices.
$\omega$ is said to be {\em horizontally homogeneous} if all the
coefficients $\omega_A\in {\A}_{hor}$ have  the same degree called
the {\em horizontal degree} and denoted by $\deg_h$.

   The component $\omega_\emptyset\in {\A}_{hor}$ of $\omega\in {\A}$ is called the {\em horizontal component}, the map  $\omega \mapsto \omega_{\emptyset}$ will be denoted by $h$ and will be named the {\em horizontal projection}.

\begin{remark}{\rm It is not difficult to see that the horizontal projection
can be explicitly described by}
 \[
   h=\prod_k\left({\bf
   1}-{\tth}^k\otimes{\I}_k\right)=\exp(-{\tth}^k\otimes{\I}_k).
\]
\label{remark: proj}
\end{remark}

 We can now define the covariant derivative of  the connection $\Theta$ as the composition
\[
\nabla =h\circ d.
\]
 A simple computation shows that
\be
\nabla {\tth}^i ={\tom}^i\;\;\;\;{\rm (Maurer-Cartan)} \label{eq:
MC} \ee and
\be
\nabla {\tom}^i =0\;\;\;\;{\rm (Bianchi)} \label{eq: Bia} \ee

 Set
\[
{\cw}: W_G \ra {\A},\;\;\;\theta^i \mapsto {\tth}^i, \;\;\Omega^i
\mapsto {\tom}^i.
\]
This map is independent of the basis $(e_i)$ and it is called {\em
the Chern-Weil correspondence}.  The following result  explains
the universal role played by the``exotic'' structure of $W_G$.
\begin{proposition}{\rm  The Chern-Weil correspondence induced by a connection is a morphism of $G$-operations.
 Moreover,  given two connections ${\tth}_i,\;i=0,1$ on ${\A}$
the corresponding Chern-Weil maps ${\cw}_i$ are homotopic as
morphisms of cochain complexes.} \label{prop: wtrans}
\end{proposition}

\noindent{\bf Proof}\hspace{.3cm}  The first part is left to the
reader.  To prove the second part we use a familiar trick from the
theory of characteristic classes.

Form the algebra $\hat{\A}=\Omega^*({\bR})\otimes {\A}$. (If
${\A}$ were the algebra of differential forms  $\Omega^*(M)$ on a
smooth manifold $M$ then $\hat{\A}\cong \Omega^*({\bR}\times M)$.)
Clearly $\hat{\A}$ is a $G$-operation  and $\hat{\theta}
=(1-t){\tth}_0+t{\tth}_1$ defines a connection on $\hat{\A}$.
Denote  by $\Psi_i$ $(i=0,1)$ the the maps $\Psi_i:\hat{\A} \ra
{\A}$ defined by the localizations at $t=1$
\[
\Omega^*({\bR})\ra {\bR}, \;\;f(t)\mapsto f(i),\;\;dt\mapsto 0.
\]
We have a fiberwise integration morphism
\[
\int_I:\hat{\A}\ra {\A}
\]
defined by
\[
\int_I f(t)\otimes\omega =\left\{
\begin{array}{rcc}
0 & {\rm if} & f\in \Omega^0({\bR}) \\
\left(\int_0^1f(t)\right)\omega & {\rm if} & f\in \Omega^1({\bR})
\end{array}
\right.
\]
The fundamental theorem of calculus implies  immediately the
following homotopy formula
\[
\forall \hat{\omega}\in \hat{\A}:
\;\;\Psi_1\hat{\omega}-\Psi_0\hat{\omega}=d\int_I\hat{\omega}
+\int_I\hat{d}\hat{\omega}
\]
where $\hat{d}$ is the exterior derivative in $\hat{\A}$ defined
by
\[
\hat{d}=dt\frac{\partial}{\partial t}\otimes {\bf 1}+\epsilon
\otimes d.
\]
$d$ is the exterior derivative in ${\A}$ while $\epsilon$ is the
s-grading operator in $\Omega^*({\bR})$. From the equalities
${\cw}_i=\Psi_i\circ \hat{\cw}$ we deduce {\em  Weil's
transgression formula}
\be
{\cw}_1-{\cw}_0=d\int_I\hat{cw}+\int_I\hat{d}\hat{\cw}=d\int_I\hat{\cw}
+\int_I\hat{\cw}d_W. \label{eq: wtransgress} \ee Thus the map
\[
K=K({\tth}_1, {\tth}_0)=\int_I \hat{\cw} :W^*_G\ra {\A}^{*-1}
\]
is a cochain homotopy connecting ${\cw}_0$ to ${\cw}_1$. $\Box$

\begin{remark}{\rm (a) It is instructive to compute $K(\theta^i)$ and
$K(\Omega^j)$. We have
\[
K(\theta^i) =\int_I\hat{\theta}^i=0.
\]
To compute $K(\Omega^j)$ we need to compute the curvature
$\hat{\Omega}$ of $\hat{\theta}$. Set
$\dot{\tth}={\tth}_1-{\tth}_0$. We have
\[
\hat{\Omega}=\hat{d}\hat{\theta}+\frac{1}{2}[\hat{\theta},\hat{\theta}]
\]
\[
=dt\otimes\dot{\tth} +d{\tth}_0 +\frac{1}{2}[{\tth}_0,{\tth}_0]
+td\dot{\tth} +
t[{\tth}_0,\dot{\tth}]+\frac{t^2}{2}[\dot{\tth},\dot{\tth}]
\]
\[
=dt\otimes \dot{\tth} +{\tom}_0  + t\left(d{\tth}_1
+[{\tth}_1,{\tth}_0] -{\tom}_0
-\frac{1}{2}[{\tth}_0,{\tth}_0]\right) +\frac{t^2}{2}[\dot{\tth},
\dot{\tth}]
\]
\[
=dt\otimes \dot{\tth} +(1-t){\tom}_0 + t\left(d{\tth}_1
+[{\tth}_1,{\tth}_0]
 -\frac{1}{2}[{\tth}_0,{\tth}_0]\right) +
\frac{t^2}{2}[\dot{\tth}, \dot{\tth}].
\]
where ${\tom}_0$ is the curvature of  ${\tth}_0$.  Set
\[
X(t)=X(t,{\tth}_1,{\tth}_0)= (1-t){\tom}_0 + t\left(d{\tth}_1
+[{\tth}_1,{\tth}_0]
 -\frac{1}{2}[{\tth}_0,{\tth}_0]\right) +
\frac{t^2}{2}[\dot{\tth}, \dot{\tth}].
\]
 Thus
\[
K(\Omega) =\int_I
\left(X(t)+dt\otimes\dot{\tth}\right)={\tth}_1-{\tth}_0.
\]
More generally if $P\in S{\gog}^*$ and $Q\in \Lambda{\gog}^*$ we
have
\[
K(P\otimes Q)=\int_IP(\hat{\Omega})\otimes Q(\hat{\theta})
\]
\[
=\int_I \left\{P( X(t)+ dt\otimes \dot{\tth})\otimes
Q(\hat{\theta})\right\}.
\]
Using the Taylor expansion of $P$ at $X(t)$ we get
\be
=\int_0^1\left(\sum_i ({\tth}_1^i-{\tth}_0^i)\frac{\partial
P}{\partial \Omega^i}(X(t))\otimes Q(\hat{\theta})\right)dt.
\label{eq: wtrans} \ee In particular, this shows that (i) $K$
commutes with the $G$ action and (ii) $K(W_{bas})\subset
{\A}_{bas}$. We call such a homotopy a {\em basic homotopy}.
 We will use the symbol ``$\simeq_b$'' to denote the basic homotopy   equivalence relation. In
 particular, we say that two $G$-operations ${\A}, {\B}$ are b-homotopic if there exist  morphisms $f:{\A}\ra {\B}$  and $g:{\B}\ra {\A}$ such that $f\circ g\simeq_b {\bf id}$ and $g\circ f \simeq_b {\bf id}$. We write this as ${\A}\simeq_b {\B}$.

(b) The above proposition shows that  we could   define the notion
of connection as a morphism of $G$-operations $W_G\ra {\A}$.  We
see that $W_G$ is extremely rigid since  for any $G$-operation the
collection $[W_G,{\A}]_b$ of classes of morphisms of
$G$-operations  modulo basic homotopies is very small. It consists
of at most one element.  A $G$-operation ${\cal W}$ equipped with
a $G$-connection satisfying the above rigidity condition (i.e.
$[{\cal W}, {\cal A}]_b$ consists of at most one element for any
$G$-operation ${\cal A}$) will be called an {\em universal
$G$-operation}.  Note also the similarity between this result and
the topological one: two continuous maps $f_1,f_2:B\ra BG$ which
induce isomorphic principal $G$-bundles  are homotopic.

(c)  The  proof of the above proposition continues to hold in the
following more general form:  any two equivariant morphisms
$\phi_i:W_G\ra {\A}$ of graded differential algebras are homotopic
as cochain maps. In particular, this shows that $(W_G,d_W)$ is
acyclic i.e. $H^k(W_G)=0$ for $k>0$.} \label{remark: wtransgress}
\end{remark}

\section{The basic cohomology of a $G$-operation with connection}
\setcounter{equation}{0}

As we explained in the previous subsection,  the  $G$-operations
with connections  represent the algebraic counterpart of a smooth
manifold  $M$ on which $G$ acts freely. In such a case, the
(Borel) equivariant cohomology  of $M$ is naturally isomorphic
with the ordinary cohomology of the quotient
\[
H^*_G(M)\cong H^*(M/G).
\]
In the subsection we will establish the algebraic counterpart of
this result. In fact, we will deal with a more general situation.

Assume we are given the following collection of  data.

\noindent $\bullet$ A Lie group $G$ and a closed normal subgroup
$N\subset G$. Set $Q=G/N$. Since $N$ is invariant under the
adjoint of $G$ there is an induced action on ${\gn}$-the Lie
algebra of $N$ and in particular, ${\gn}$ is a Lie algebra ideal
of ${\gog}$.

\noindent $\bullet$ A $G$-operation $({\A}, D, {\L}, {\I})$
equipped with a $G$-invariant {\it $N$-connection} i.e. a
$G$-equivariant morphism of {\it $N$-operations}
\[
{\tth}:W_N \ra {\A}.
\]
By regarding ${\A}$ as an $N$-operation we can form the subalgebra
\[
{\B}={\A}_{N,bas}=\{\omega\in {\A}\; ;\; \omega\; \;{\rm is}\; N\;
{\rm invariant},\;{\I}_X\omega =0,\;\forall X\in {\gn}\}.
\]
The $G$-operation structure on ${\A}$ induces a residual
$G/N$-operation structure  on ${\B}$. Note that we have an
inclusion
\[
j:W_Q\otimes {\B}\ra W_G\otimes {\A}
\]
such that
\[
j(W_Q\otimes {\B})_{Q,bas} \subset (W_G\otimes {\A})_{G,bas}.
\]
The geometric intuition behind this algebraic situation is that of
a   smooth $G$-space $E$ such that the action of $N$ is free.  In
Borel cohomology we have an isomorphism
\[
H^*_G(E)\cong H^*_Q(E/N).
\]
We will establish the algebraic analogue of this result.
\begin{theorem} {\bf (Cartan)} {\rm  The inclusion $j$ induces an isomorphism}
\[
WH_Q^*({\B})\cong WH^*_G({\A}).
\]
\label{th: Cartan}
\end{theorem}

\noindent {\bf Proof}\hspace{.3cm} Our proof will be a simple
application of Weil's transgression trick. For different
approaches we refer to \cite{GHV}, Chap. VIII, \cite{Du},
\cite{DV} or \cite{KV}.

We will construct a $G$-connection on $W_G\otimes{\A}$ starting
from the $G$-equivariant $N$-connection ${\tth}\in  {\A}^1\otimes
{\gn}$.

Define the linear map
\[
\mu:{\gog}\ra {\A}^0\otimes {\gn},\;\;\mu(X)=-{\I}_X{\tth}.
\]
$\mu$ is called the {\em moment map} of the $G$-equivariant
connection ${\tth}$. We can also regard it as  an element of
${\gog}^*\otimes {\A}^0\otimes {\gn}$.

\begin{lemma}{\rm The moment map $\mu : {\gog}\ra {\A}^0\otimes {\gn}$ is
$G$-equivariant.} \label{lemma: mu}
\end{lemma}

\noindent{\bf Proof}\hspace{.3cm} Regard ${\tth}$ as a
$G$-equivariant map
\[
{\gog}^*\ra {\gn}^*\ra {\A}^1.
\]
For each $X\in {\gog}$ regard $\mu(X)$ as a map
\[
{\gog}^*\ra {\gn}^* \ra {\A}^0.
\]
The equivariance of $\mu$ is equivalent to
\[
\mu ( Ad_gX)= g\mu(X) Ad^*_{g^{-1}}
\]
where $Ad^*$ denotes the coadjoint action of $G$. We have
\[
\mu ( Ad_gX)=-{\I}_{Ad_gX}{\tth}=-g{\I}_X
g^{-1}{\tth}=-g{\I}_X{\tth}Ad^*_{g^{-1}}=g\mu(X) Ad^*_{g^{-1}}.
\]
(The second equality is the $G$-equivariance of ${\I}$ while the
third equality is the $G$-equivariance of ${\tth}$.)  $\Box$

\bigskip

Define ${\bf q}:{\gog}\ra {\A}^0\otimes {\gn}$ as
\[
{\bf q}(X)=X+\mu(X).
\]
Note that ${\bf q}(X)=0$ for $X\in {\gn}$ so that {\bf q} descends
to a map
\[
{\bf q}={\goq}={\gog}/{\gn} \ra {\A}^0\otimes {\gn}.
\]
Set $\Xi= {\bf q} +{\tth}$. Note that
\[
{\bf q}\in {\gog}^*\otimes {\A}^0 \otimes {\gn} \subset
W^1_G\otimes {\A}^0\otimes {\gog}
\]
and
\[
{\tth} \in {\A}^1\otimes {\gn} \subset {\A}^1 \otimes{\gog}.
\]
Thus $\Xi \in (W_G\otimes {\A})^1$.

\begin{lemma}{\rm  $\Xi$ defines a $G$-connection on $W_G\otimes {\A}$}
\end{lemma}

\noindent{\bf Proof}\hspace{.3cm} For $X\in {\gog}$ denote by
$I_X$ the contraction by $X$ in $W_G\otimes {\A}$. Then
\[
I_X\Xi= {\bf q}(X) +{\I}_X{\tth} =X.
\]
The $G$-invariance of $\Xi$ now follows from the $G$-invariance of
${\bf q}$ and ${\tth}$.  $\Box$

\bigskip

The $G$-operation $W_G\otimes {\A}$ admits the tautological
connection
\[
{\bf 1}:W_G \mapsto W_G \otimes {\A},\;\;w\mapsto w\otimes 1.
\]
Denote by $K=K({\bf 1}, \Xi)$ the Weil transgression
\[
K:W_G \ra W_G\otimes{\A}
\]
so that for all $w\in W_G$
\[
w-\Xi w=\delta K w + K d_W w
\]
where $\delta$ denotes the  exterior derivation in $W_G\otimes
{\A}$. Now define
\be
T_0:W_G\otimes {\A} \ra W_G\otimes {\A},\;\;w\otimes a=(\Xi
w)\cdot  a \label{eq: red} \ee
\[
T_1={\bf id}: W_G\otimes {\A}\ra W_G\otimes {\A}
\]
and
\[
{\cal K}: W_G\otimes {\A}\ra W_G\otimes {\A},\;\;w\otimes a
\mapsto Kw\cdot a.
\]
Then for all $x\in W_G\otimes {\A}$
\[
x-T_0x=T_1x-T_0x=\delta{\cal K}x +{\cal K}\delta x.
\]
Both $T_0$ and $T_1$ are morphisms of $G$-operations and ${\cal
K}$ is a {\em basic homotopy}. Also note that
\[
T_0(W_G\otimes {\A})\subset W_Q\otimes {\A}
\]
and
\[
T_0(W_G\otimes {\A})_{G,bas} \subset (W_Q\otimes{\B})_{Q,bas}.
\]
Moreover, along the basic subalgebras $T_0\circ j={\bf id}$. In
(the basic) cohomology $T_0$ is bijective since it is homotopic to
the identity. This completes the proof of Theorem \ref{th:
Cartan}.   $\Box$

\begin{remark}{\rm    The reduction theorem we have just proved generalizes as follows. Consider a $G$-operation  ${\W}$. Then  the transgression trick in the proof of Theorem \ref{th: Cartan}  can be used to show the statements below are equivalent.

\noindent (i) ${\W}$ is universal (in the sense defined in  Remark
\ref{remark: wtransgress} (b)).

\noindent (ii) Any morphism of $G$-operations ${\vfi}:{\W}\ra
{\A}$  induces a b-homotopy equivalence
\[
{\W}\otimes {\A}\simeq_b {\A}.
\]

Note in particular that if ${\W}_0$, ${\W}_1$ are two universal
$G$-operations and $\alpha:{\W}_0\ra {\W}_1$ is a morphism  then
any morphism  $\tau: {\W}_1\ra {\W}_0\otimes{\B}$ induces  a
b-homotopy equivalence
\be
{\W}_0\otimes {\B} \simeq_b {\W}_1 \otimes {\B}. \label{eq: bhomo}
\ee Indeed,   by (ii) $\tau$ induces a b-homotopy equivalence
\[
{\W}_0\otimes {\W}_1\otimes {\B}\simeq_b {\W}_1\otimes
\left({\W}_0\otimes {\B}\right)\simeq_b  {\W}_0\otimes {\B}
\]
(since ${\W}_1$ is universal) and on the other hand $\alpha :
{\W}_0\ra {\W}_1 \otimes {\B}$ induces the b-homotopy equivalence
\[
{\W}_0\otimes {\W}_1\otimes  {\B} \simeq_b {\W}_1\otimes {\B}
\]
since ${\W}_0$ is universal. If we take   ${\W}_0=W_{G/N}$ and
${\W}_1=W_G$  (note that these Weil algebras are clearly
universal $G$-operations) then the  equivalence (\ref{eq: bhomo})
 is precisely the content of Theorem \ref{th: Cartan}.}
 \label{rem: proof}
\end{remark}

\begin{corollary}{\rm Let $E$ be a smooth  $G$-space. If $N$ acts freely
on $E$ then}
\[
WH^*_G(E)\cong H_Q^*(E/N).
\]
\end{corollary}

It is instructive to describe the reduction isomorphism
\[
WH^*_G(E)\stackrel{\cong}{\ra} WH^*_Q(E/N)
\]
using the Cartan model.  Denote by $W_{\tth}$ the Weil model
description of the above reduction isomorphism (defined in
(\ref{eq: red})). We denote by $C_{\tth}$ its correspondent in the
Cartan model.   Denote by $\phi_G$ (resp. $\phi_Q$) the Kalkman
isomorphism (cf. (\ref{eq: kalk}))
\[
WH^*_G\ra CH^*_G\;\;\;({\rm resp.} WH^*_Q\ra CH^*_Q).
\]
We then  have
\[
C_{\tth}=\phi_Q\circ W_{\tth}\circ \phi^{-1}_G.
\]
To get a better feeling on the  structure of $C_{\tth}$ we will
work in local coordinates. Choose a basis $(e_i)$ of ${\gn}$ and
then extend it to a basis $ \{e_i\, ;\, f_a\} $ of ${\gog}$. Via
the natural projection ${\gog}\ra {\gog}/{\gn}$ the collection
$(f_a)$ induces a basis of  ${\gog}/{\gn}$ which we continue to
denote by the same symbols.  Denote the dual basis of  $ \{e_i\,
;\, f_a\} $ by   $\{ \theta^i \, ; \, {\vfi}^a\}$. We can regard
$(\theta^i)$ as a basis of  ${\gn}^*$ and $({\vfi}^a)$ as a basis
of ${\goq}^*$.   We denote the image of $\theta^i$ in ${\gn}^*$ by
$\Omega^i$ and the  image of  ${\vfi}^a$ in $S{\goq}^*$ by
$\Psi^a$. Set $\theta=\theta^i\otimes e_i$,
${\vfi}={\vfi}^a\otimes f_a$, $\Omega= \Omega^i\otimes e_i$ and
$\Psi=\Psi^a\otimes f_a$.Then  any element in $(S{\gog}^*\otimes
\Omega^*(E))^G$ is a $G$-equivariant polynomial map $P:{\gog}\ra
\Omega^*(E)$ which we will schematically describe  it as
$P=P(\Omega \oplus \Psi)$. Then
\[
\phi_G^{-1}P(\Omega\oplus \Psi)=\exp(-{\vfi}^a\otimes {\I}_a)\exp
(-\theta^i\otimes {\I}_i)P(\Omega\oplus \Psi).
\]
The map $\Xi: W_G\ra W_G\otimes\Omega^*(E)$ is determined from the
assignments
\[
\theta^i\mapsto {\tth}^i-{\tth}^i(f_a){\vfi}^a,\;\;{\vfi}^a\mapsto
{\vfi}^a.
\]
  If we define $\Xi(e_i)=e_i$ and $\Xi(f_a)=f_a$ then we can rewrite
\be
\Xi(\theta\oplus {\vfi})=
{\tth}+{\vfi}-{\tth}(f_a){\vfi}^a={\tth}+{\vfi}-{\tth}\circ{\vfi}.
\label{eq: reduc} \ee Moreover
\[
\Xi(\Omega \oplus \Psi) =\delta\Xi(\theta\oplus {\vfi})
+\frac{1}{2}\left[ \Xi(\theta\oplus
{\vfi}),\Xi(\theta\oplus{\vfi})\right]
\]
\[
=d{\tth}+d_W{\vfi}
+\frac{1}{2}\left\{[{\tth}+{\vfi},{\tth}+{\vfi}] +
[{\tth}\circ{\vfi},{\tth}\circ{\vfi}] \right\}-
[{\tth}+{\vfi},{\tth}\circ{\vfi}] -\delta({\tth}\circ {\vfi})
\]
\[
=d{\tth}+\frac{1}{2}[{\tth},{\tth}]+d_W{\vfi}+\frac{1}{2}[{\vfi},{\vfi}]
+ [{\tth},{\vfi}] +
\frac{1}{2}[{\tth}\circ{\vfi},{\tth}\circ{\vfi}]
-[{\tth}+{\vfi},{\tth}\circ{\vfi}]-\delta({\tth}\circ {\vfi})
\]
\be
={\tom} +\Psi  + [{\tth},{\vfi}] +
\frac{1}{2}[{\tth}\circ{\vfi},{\tth}\circ{\vfi}]
-[{\tth}+{\vfi},{\tth}\circ{\vfi}]-\delta({\tth}\circ {\vfi}).
\label{eq: reduc1} \ee On the other hand
\be
W_{\tth}\circ \phi^{-1}_GP(\Omega \oplus \Psi)
=\exp(-{\vfi}^a\otimes {\I}_a)\exp(-\Xi(\theta^j)\otimes
{\I}_j)P(\Xi(\Omega\oplus \Psi)). \label{eq: reduc2} \ee Since
$W_{\tth}\circ\phi^{-1}_GP$ is a $Q$-basic element of $W_Q\otimes
\Omega^*(E/N)$ the action of $\phi_Q$ on this element is
determined according to Remark \ref{remark: brst} by setting
${\vfi}=0$ in (\ref{eq: reduc2}). Using  the equalities (\ref{eq:
reduc}) and (\ref{eq: reduc1}) we deduce
\be
C_{\tth}P=\phi_Q\circ W_{\tth}\circ \phi_G^{-1}P=
\exp(-{\tth}^j\otimes {\I}_j)P({\tom}+\Psi-  \delta({\tth}\circ
{\vfi})\!\mid_{{\vfi}=0}  \,). \label{eq: reduc3} \ee We need to
understand the term $ \delta({\tth}\circ {\vfi})\!\mid_{{\vfi}=0}
$.   We have
\[
\delta ({\tth}^i(f_a){\vfi}^a)=d({\tth}^i(f_a)){\vfi}^a
+{\tth}^id_W{\vfi}^a
\]
\[
=d({\tth}^i(f_a)){\vfi}^a  +{\tth}^i(f_a)\{\Psi^a -{\cal
Q}({\vfi})\}
\]
where ${\cal Q}$ denotes a quadratic term in the ${\vfi}$'s. Thus
when setting ${\vfi}=0$ we get
\[
\delta({\tth}^i(f_a){\vfi}^a\circ {\vfi})\!\mid_{{\vfi}=0}
={\tth}^i({f_a})\Psi^a.
\]
Symbolically
\[
\delta({\tth}\circ {\vfi})\!\mid_{{\vfi}=0}
={\tth}(\Psi)=-\mu(\Psi)\in (S^1{\goq}^*\otimes
\Omega^0(E))\otimes {\gn}.
\]
(Recall that $\mu$ denotes the moment map of the connection.)
Note that the differential form components in $\mu(\Psi)$ are
$N$-basic so they can be regarded as forms on the basis $E/N$.
Substituting this back in (\ref{eq: reduc3}) we get
\[
C_{\tth}P=\exp(-{\tth}\otimes {\I}_j)P(\Psi+\mu(\Psi)+{\tom}).
\]
The exponential factor is precisely the  horizontal projection
$h_{\tth}:\Omega^*(E)\ra \Omega^*(E/N)$ defined by the
$N$-connection ${\tth}$. On the other hand  the term
${\tom}+\mu(\Psi)\in (\Omega^2(E/N)\oplus (S^1{\goq}^*\otimes
\Omega^0(E/N))\otimes {\gn}$ is already $Q$-basic. It is called
the {\em equivariant curvature} of the connection ${\tth}$ and
will be denoted by ${\tom}_Q$. Note that ${\tom}_Q$ is an element
of degree  2 in ${\tom}_Q\in (S{\goq}^*\otimes
\Omega^*(E/N))^Q\otimes {\gn}$. Thus
\be
(C_{\tth}P)(\Psi)=h_{\tth}P(\Psi+{\tom}_Q). \label{eq: reduc4} \ee
We still need to give an accurate definition of the
right-hand-side term above. For any $X \in {\gog}$ define $P(X
+{\tom}_Q)$ imitating the  Taylor expansion at $X$
\[
P(X+{\tom}_Q)=\exp({\tom}_Q^i\partial_i)P(X)
\]
where ${\tom}_Q={\tom}^i_Q\otimes e_i \in (S{\goq}^*\otimes
\Omega^*(E))^2\otimes {\gn}$ while $\partial_i$ denotes the
partial derivative along the direction $e_i\in {\gn}\subset
{\gog}$. Note that $P(X+{\tom}_Q)=P(X+\mu(X)+{\tom})=P({\tom})$
for all $X\in {\gn}$ so that $P(X+{\tom}_Q)$ descends to a well
defined map ${\goq}\ra \Omega^*(E/N)$.  Thus the polynomial in the
right-hand-side of (\ref{eq: reduc4}) should be rather viewed as
an $\Omega^*(E/N)$-valued polynomial on ${\gog}$ which descends to
a polynomial  on ${\goq}={\gog}/{\gn}$.

In particular, to any $G$-invariant polynomial $P\in S{\gog}^*$
one can associate an  equivariantly closed element
\[
\Psi \mapsto P(\Psi+{\tom}_Q)\in (S{\goq}^*\otimes
\Omega^*(E/N))^Q.
\]
 This element clearly depends on the connection but  its
equivariant cohomology class does not.  It will be  denoted by
$P(E)\in CH^*_Q((E/N)$ and will be called the {\em equivariant
characteristic class} of $E\ra E/N$ corresponding to $P$.  Note
that when $G=N$  this correspondence is none other than the
traditional Chern-Weil construction of the characteristic classes
of a  principal $G$-bundle.

\end{document}